%% file: CameraReady_LZeune.tex
\documentclass[runningheads,a4paper]{llncs}

\usepackage{amssymb}
\usepackage{graphicx}
\usepackage{xcolor}
\usepackage{bbold}
\usepackage{gensymb}
\usepackage{amsmath}
\usepackage{lipsum}
\usepackage{amsfonts}
\usepackage{epstopdf}
\usepackage{authblk}
\usepackage{algorithm}
\usepackage{algpseudocode}
\usepackage{subfig}
\usepackage{bbm}
\input{commands}%

\usepackage{url}
\urldef{\mailsa}\path|{l.l.zeune, s.a.vangils, c.brune}@utwente.nl|
\urldef{\mailsb}\path|l.w.m.m.terstappen@utwente.nl|
\newcommand{\keywords}[1]{\par\addvspace\baselineskip
\noindent\keywordname\enspace\ignorespaces#1}

\begin{document}

\mainmatter  % start of an individual contribution

% first the title is needed
\title{Combining Contrast Invariant L1 Data Fidelities with Nonlinear Spectral Image Decomposition}

% a short form should be given in case it is too long for the running head
\titlerunning{Contrast Invariant L1 Data Fidelities for Spectral Analysis}

% the name(s) of the author(s) follow(s) next
%
% NB: Chinese authors should write their first names(s) in front of
% their surnames. This ensures that the names appear correctly in
% the running heads and the author index.
%
\author{Leonie Zeune\inst{1,}\inst{2}\and Stephan A van Gils\inst{1}\and Leon WMM Terstappen\inst{2}\and Christoph Brune\inst{1}}
\authorrunning{L. Zeune, S.A. van Gils, L.W.M.M. Terstappen and C. Brune}
% (feature abused for this document to repeat the title also on left hand pages)

% the affiliations are given next; don't give your e-mail address
% unless you accept that it will be published
\institute{Department of Applied Mathematics\\
\mailsa
\and
Department of Medical Cell BioPhysics\\
\mailsb \\
University of Twente, 7522 NB Enschede, The Netherlands\\
}

%
% NB: a more complex sample for affiliations and the mapping to the
% corresponding authors can be found in the file "llncs.dem"
% (search for the string "\mainmatter" where a contribution starts).
% "llncs.dem" accompanies the document class "llncs.cls".
%

%\toctitle{Lecture Notes in Computer Science}
%\tocauthor{Authors' Instructions}
\maketitle

\begin{abstract}%
\vspace*{-1em}%
This paper focuses on multi-scale approaches for variational methods and corresponding gradient flows. Recently, for convex regularization functionals such as total variation, new theory and algorithms for nonlinear eigenvalue problems via nonlinear spectral decompositions have been developed. Those methods open new directions for advanced image filtering. However, for an effective use in image segmentation and shape decomposition, a clear interpretation of the spectral response regarding size and intensity scales is needed but lacking in current approaches. In this context, $L^1$ data fidelities are particularly helpful due to their interesting multi-scale properties such as contrast invariance. Hence, the novelty of this work is the combination of $L^1$-based multi-scale methods with nonlinear spectral decompositions. We compare $L^1$ with $L^2$ scale-space methods in view of spectral image representation and decomposition. We show that the contrast invariant multi-scale behavior of $L^1-TV$ promotes sparsity in the spectral response providing more informative decompositions. We provide a numerical method and analyze synthetic and biomedical images at which decomposition leads to improved segmentation.
\end{abstract}
\section{Introduction}
%%%%%%%%%%%%%%%%%%%%
\vspace{-.5em}%
In imaging science, the solution of inverse problems is often addressed by the modeling and analysis of variational methods of the form:
\begin{equation}\label{eq:varMeth}%
	\min_{u} ~ \frac{1}{p} \norm{u-f}_{L^p}^{p} + \alpha J(u)
\end{equation}
where $f$ denotes a noisy signal, $u$ is a desired image function defined on $\Omega$ in $\R^2$. The data fidelity measures the residual in the $L^p$ norm and $J$ is a regularization functional that has a weighting parameter $\alpha \geq 0$. In this paper, we concentrate on image denoising methods with convex and one-homogeneous regularization functionals $J$ which can address image decomposition and segmentation adequately. More specifically, we focus on total variation regularization $J(u) = TV(u)$ and norms with $p=1$ versus $p=2$. One option to generate a scale-space method is a gradient flow, based on the functional \eqref{eq:varMeth} with initial condition $u(0,x) = f(x)$ and subdifferential inclusions forming the doubly nonlinear evolution equation:
\begin{equation}\label{eq:tvflow}%
	\vspace{-.5em}%
	0 \in \partial \norm{\partial_t u(t,x)}_{L^p} +  \partial J(u(t,x)) .
\end{equation}%
In the case of $p=2$, this simplifies to well-known gradient flows of the form:
\begin{equation}\label{eq:tvflowL2}%
	\vspace{-.5em}%
	\partial_t u = - q(u) \quad \text{with } q \in \partial J(u) .
\end{equation}%
The regularization parameter $\alpha$ is now hidden in the time dependency. Corresponding inverse scale space gradient flows can be constructed via Bregman distances \cite{Osher2005,Brune2011}. The analysis of linear eigenvalue problems and spectral decompositions, e.g. via the Fourier transform, is a well-known and widely used theory in the fields of signal, image and graph-based data processing. Due to the continuing success of nonlinear regularization functionals in imaging, there is strong interest in generalizing spectral theory to the nonlinear case. The general idea is to examine solutions to the nonlinear eigenvalue problem:
\vspace*{-.5em}%
\begin{equation}\label{eq:eigenval}
	\lambda u \in \partial J(u) .
\end{equation}%
In \cite{Gilboa2013,Gilboa2014}, Gilboa introduced the idea of nonlinear spectral decompositions for the TV transform. By transferring solutions (eigenfunctions) of \eqref{eq:eigenval} to sparse peaks in a spectral domain, the idea of advanced filters, suppressing or enhancing image components similar to the Fourier transform, came about. This concept was studied for one-homogenous functionals \cite{Burger2015,Burger2016,Gilboa2015} and scale-space flows of the form \eqref{eq:tvflowL2} with $p=2$. In this way, a decomposition of the input signal $f$ into significant components is possible while an exact representation of $f$ can still be guaranteed.

\begin{figure}[htb]
	\centering%
	\begin{minipage}{.13\textwidth}
		\vspace{-0.21\textheight}
		\includegraphics[width=\textwidth]{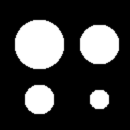}%
	\end{minipage}%
	\quad~%
  %\subfloat[]{
	\includegraphics[width=0.37\textwidth]{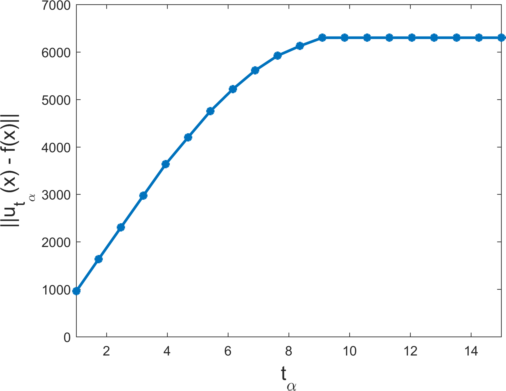}
	%}
	\quad%
  %\subfloat[]{
	\includegraphics[width=0.37\textwidth]{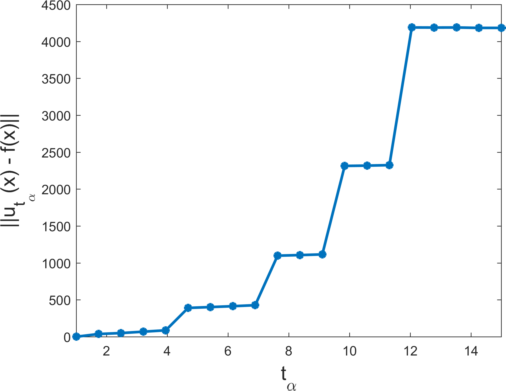}
	%}
	\caption{\textit{Scale-spaces.} (left) input image $f$, (middle) $L^2-TV$, (right) $L^1-TV$.}
	\label{fig:scaleSpace}%
	%\vspace*{-2em}%
\end{figure}

For $p=1$, the gradient flow in \eqref{eq:tvflow} is interesting and more challenging than for $p=2$. Another way to obtain a forward scale-space is to construct a sequence of variational problems of the form \eqref{eq:varMeth} with increasing regularization parameter $t_\alpha$ replacing the fixed $\alpha$. Here, the scale parameter $t_\alpha$ corresponds to the time variable $t$ used in \eqref{eq:tvflow}. In this paper, we will focus on this type of model particularly in view of $p=1$. From pioneering works on nonlinear $L^1$ filtering \cite{Nikolova2002,Chan2005,Aujol2006,Yin2007}, it is well known that such variational reconstruction methods share interesting multi-scale properties including contrast invariance. For this reason, $L^1$ data fidelities have been successfully used for advanced image reconstruction techniques \cite{Dong2009,Wu2011,Yuan2011}, vector field estimation \cite{Zach2007} and image decompositions regarding texture \cite{Aujol2006,Haddad2007}. The special multi-scale behavior becomes clear in Figure \ref{fig:scaleSpace}. In the $L^1-TV$ case shown on the right, the contrast invariance leads to plateaus in the scale-space graph indicating an abrupt disappearance of TV-eigenshapes regarding the input image on the left.\\
Motivated by \cite{Chan2005,Yin2007,Duval_2009} and nonlinear spectral methods \cite{Gilboa2015,Zeune2016}, the main goal of this work is to study $L^1$ versus $L^2$ in view of the sparsity of nonlinear spectral decompositions. Does $L^1$ imply sparsity and hence a more \textit{informative spectral response}? Can we expect more \textit{reliable image decompositions} via backtransformation facilitating improved image segmentation? What happens if complex shapes, e.g. nonconvex or \textit{compositions of eigenfunctions}, are involved?
%
%%%%%%%%%%%%%%%%%%%%

%\vspace*{-2em}%
\section{Modeling}
 \vspace*{-0.5em}
%%%%%%%%%%%%%%%%%%%%
In the following section we first give a short overview of the eigenshapes of the total variation functional and calibrable sets before we introduce the spectral framework for $L^2-TV$ denoising in more detail. Finally, we show why a combination of the spectral framework with $L^1-TV$ denoising seems promising and how the spectral framework can be adapted in this case.
\vspace*{-0.5em}%
\subsection{Geometry: Eigenfunctions and Calibrable Sets}
In the introduction in Figure \ref{fig:scaleSpace} we have already seen an input image composed of TV eigenfunctions. A  general geometric description of TV eigenfunctions is given in \cite{Bellettini2002,Caselles2007,Chambolle2016} in terms of calibrable sets, more precisely by convex sets which are Cheeger in themself. An indicator function $\chi_C(x)$ of a convex and connected set $C$ with finite perimeter $Per(C)$, for which $C$ admits
\begin{equation*}\label{eq:l2eig}
	\qquad\qquad \esssup_{x \in \partial C} \kappa_C (x) \leq \frac{Per(C)}{|C|} \qquad \text{($C$ is Cheeger in itself)}
\end{equation*}
where $|C|$ denotes the area and $\kappa_C$ the curvature of $\partial C \in C^{1,1}$, is a solution of \eqref{eq:eigenval} with unit norm and therefore an eigenfunction of TV. With this geometric interpretation of eigenfunctions for TV as Cheeger sets, the role of perimeter and volume is significant for contrast invariant image decompositions. In the convex case, Duval and collaborators \cite{Duval_2009} proved that exact solutions of the $L^1-TV$ problem are given by an morphological opening followed by a simple test over the perimeter-area ratio. This fact was first published but not proved in \cite{Darbon2005}. For more complex shapes, formed by compositions of eigenfunctions, a better understanding of the $L^1$ scale-space flows is therefore a very promising direction.

\vspace*{-1em}%
\subsection{$L^{2}-TV$ based Spectral Analysis}
The nonlinear spectral analysis framework was first introduced by Gilboa in \cite{Gilboa2013,Gilboa2014} for the total variation regularization functional. A forward scale-space was constructed via the TV flow \eqref{eq:tvflowL2}. Here, the first iterate is the original data $f$ which is then smoothed in every time-step such that increasingly fine scales are removed. This concept was later generalized to more general one-homogenous regularization functionals $J(u)$ in \cite{Burger2015,Burger2016,Gilboa2015}. Moreover, it was shown that the forward scale-space can also be constructed via a variational approach by iteratively solving the ROF model \eqref{eq:varMeth} for $p=2$ and $J(u) = TV(u)$ using increasing regularization parameters $t_{\alpha}$. The idea of the spectral filtering approach is to transform the signal $f$ into a spectral domain, where both filtering of certain scales and decomposition of $f$ into significant signals is possible. Thus, eigenfunctions will be mapped onto a peak in the spectral domain. The shape of the eigenfunction is determined by the chosen regularization functional, for $J(u) = TV(u)$ the most prominent eigenfunction is a disc with radius $r$ surrounding the origin. In \cite{Strong2003}, Strong and Chan analyzed how the solution of the ROF model behaves for increasing $t_{\alpha}$ if $f$ is an eigenfunction. Thus, let $f(x) = c \cdot \mathbbm{1}_{B_{r}(0)}(x)$ be a disc of constant height $c$ and radius $r$ surrounding the origin and with a background of zero. Then the solution of the $L^2-TV$ model is given as: 
\begin{equation}
u(t_{\alpha},x) = \begin{cases} \left(c - \frac{2}{r} \ t_{\alpha}\right) \cdot \mathbbm{1}_{B_{r}(0)}(x) &\mbox{if } 0 \leq t_{\alpha} < \frac{cr}{2} \\
0 &\mbox{otherwise.}
\end{cases}
\end{equation}
That means that even in the noise-free case, the reconstructed solution $u(t_{\alpha},x)$ never reaches the true value $f$ and for increasing regularization parameter $t_{\alpha}$  the disc height decreases. Solutions are corrupted by a systematic contrast loss that is dependent on the regularization strength but also on the radius $r$ and the height $c$ of input data $f$. To transform such eigenfunctions to a peak in the spectral domain, Gilboa defined the spectral transform function $\phi(t,x)$ and the spectral response function $S(t)$ as:
\begin{equation}\label{eq:tvtransform}
\phi (t,x) = u_{tt}(t,x)\cdot t \quad\text{and}\quad S(t) = \norm{\phi(t,x)}_{L^1}.
\end{equation}
The definition of $\phi$ allows, under certain conditions, that the original signal $f$ can be reconstructed via: 

\begin{equation*}
f(x) = \int_0^{\infty} \phi(t,x) dt + \bar{f}
\end{equation*}
where $\bar{f}$ is the mean of $f$. Filtered versions of $f$ can be constructed by applying:

\begin{equation*}
f_{H}(x) = \int_0^{\infty} H(t)\phi(t,x) dt + H(\infty)\bar{f}
\end{equation*}
where $H(t)$ is the filter function. 

However, a disadvantage of the $L^2$ based spectral framework is ambiguity with respect to size and intensity scales. The method is not able to clearly differentiate size and intensity scales since the timepoint $t_d$ at which a disc disappears and a peak occurs is determined by both values together. 

\vspace*{-0.5em}%
\subsection{$L^{1}-TV$ based Spectral Analysis}
In the following section we want to combine the idea of nonlinear spectral TV analysis with the $L^1$ denoising model: 
\begin{equation}\label{eq:l1tv}%
	\min_{u} ~  \norm{u-f}_{L^1} + t_{\alpha} TV(u).
\end{equation}
As mentioned earlier, this model shows very interesting multi-scale decomposition behavior and seems therefore very suitable to be combined with the spectral approach. In \cite{Chan2005}, Chan and Esedoglu showed that the behavior of solutions of the $L^1-TV $ model if $f$ is the TV eigenfunction $f(x) = c \cdot \mathbbm{1}_{B_{r}(0)}(x)$ is significantly different from the $L^2$ case. Here, the solution is given by:
\begin{equation}\label{eq:l1sol}
u(t_{\alpha},x) = \begin{cases} f &\mbox{if } 0 \leq t_{\alpha} < \frac{r}{2} \\
c' \cdot f \text{ with } c'\in [0,1] &\mbox{if } t_{\alpha} = \frac{r}{2}\\
0 &\mbox{otherwise .}
\end{cases}
\end{equation}
An interesting observation is that the solution $u(t_{\alpha},x)$ is \textbf{not} dependent on the height $c$ of the disc but only on the radius $r$. That means that the $L^1-TV$ denoising approach is contrast-invariant and therefore highly suitable for decomposing a signal $f$ based on size scales. Note that the solution of \eqref{eq:l1tv} is not unique; for $f$ defined as above and $t_{\alpha} = \frac{r}{2}$ there exist, for example, an infinite number of solutions. 

Since the gradient flow \eqref{eq:tvflow} is more challenging for $p=1$ rather than $p=2$ (due to non-smoothness of the data-term), we construct a forward scale-space by taking a variational approach; in other words, we solve \eqref{eq:l1tv} for increasing regularization parameters $t_{\alpha_1} < t_{\alpha_2} < ... < t_{\alpha_N}$. To transform an eigenfunction of  $TV(u)$ onto a peak in the spectral domain, a suitable spectral transformation function is now given by:
\begin{equation}\label{eq:tvtransform}
\phi (t,x) = -u_{t}(t,x) \text{ ~fulfilling~ } f(x) = \int_0^{\infty} \phi(t,x) dt + \hat{c}.
\end{equation}
Here, the first time derivative of $u(t,x)$ is in a distributional sense and the constant $\hat{c}$ is the median of $f$. Since the disc's height does not decrease in every time step but remains the constant $c$ until it immediately decreases to 0, already the first time-derivative leads to a delta peak in the spectral transform function defined via:
\begin{equation}
S^2(t) = \langle\phi(t,x),f(x)\rangle.
\end{equation} 
This definition was first introduced by Burger and collaborators in \cite{Burger2016} and leads to an analogue to Parseval's identity:
\begin{equation}\label{eq:parseval}
\norm{f}^2 = \langle f,f\rangle = \int_0^{\infty}\langle \phi(t,x),f(x)\rangle \mathrm{d}t = \int_0^{\infty} S^2(t) \mathrm{d}t.
\end{equation}
The signal can be filtered based on size via:
\begin{equation} \label{eq:filter}
f_{H}(x) = \int_0^{\infty} H(t)\phi(t,x) dt + H(\infty)\hat{c}
\end{equation}
where $H(t)$ is again the filter function. A segmentation of objects in a certain size range (nearly independent of their intensity; only objects with intensity $\hat{c}$ cannot be found) is given by:
\begin{equation}\label{eq:seg}
f_{\text{seg},H}(x) = \left( \int_0^{\infty} H(t)\phi(t,x) dt > 0\right).
\end{equation}
Applications where this intensity-independent segmentation approach is especially helpful are described in section \ref{sec:res}.
%%%%%%%%%%%%%%%%%%%%

\vspace*{-1em}%
\section{Numerical Approach}
%%%%%%%%%%%%%%%%%%%%
\vspace*{-.5em}%
The following section introduces the numerical realization of the spectral $L^1-TV$ framework presented previously. The main component is the numerical solution of the denoising problem. To find solutions of this purely primal nonlinear minimization problem \eqref{eq:l1tv}, we make use of the first order primal dual algorithm proposed by Chambolle and Pock \cite{Chambolle2011}. As the name already suggest, the minimization scheme works with a primal dual version of \eqref{eq:l1tv} given by: 
\vspace*{-.5em}%
\eqn{\label{eq:pimduall1tv}\langle \grad u,g\rangle + \norm{u  - f}_{L^1} - t_{\alpha} \delta_{P}(g) \longrightarrow \min_{u}\max_{g}}
where $P = \left\{g: \| g\|_{\infty} \leq 1\right\}$ and $\delta_{P}(g)$ equals $0$ if $g \in P$, and equals $\infty$ if $g\notin P$.
We define $K(u) = \grad u$, $F(u) =  \norm{u }_{L^1}$ and $G(u) = \norm{u  - f}_{L^1}$.
\vspace*{-1.5em}%
\begin{algorithm}
	\caption{First-order primal-dual algorithm to solve \eqref{eq:l1tv}.}
	\label{alg:cpforl1tv}
	{\fontsize{8}{8}\selectfont % define a different font size here
	\begin{algorithmic}%
		\State \textbf{Parameters:} data $f$, reg. param. $0<t_{\alpha_1} < t_{\alpha_2} < ... < 	
		t_{\alpha_N}$, $\tau, \sigma > 0$, \\\quad \quad \quad\quad \quad \quad\ $\theta \in [0,1]$, $maxIts
		 \in \mathbb{N}$
		\State \textbf{Initialization:} $n = 0, \ u^0=0, \ p_0:=0, \bar{u}^0 = u^0$\\
		\State \textbf{Iteration:}\\
		\State \textbf{for } \big($i = 1:N$\big) \textbf{ do }
		\begin{enumerate}\itemsep5pt
			\item Set $\alpha=t_{\alpha_i}$.
		\end{enumerate}		
		\begin{quote}
		 \State  \textbf{while } \big($n<maxIts$\big) \textbf{ do }\\
		\begin{enumerate}\renewcommand{\labelenumi}{\alph{enumi})}\itemsep5pt
			\item $g^{n+1} = \text{Proj}_{\left\{\left\{g: \| g\|_{\infty} \leq 1\right\}\right\}}\left( g^n + \sigma 
			\grad \bar{u}^n \right)$.
			\item $\text{arg}_u = u^{n} + \tau\div g$.
			\item $u^{n+1}(x) = \begin{cases}
			\text{arg}_u(x) - \frac{\tau}{\alpha} &\mbox{if }\quad \text{arg}_u(x)-f(x) > 
			\frac{\tau}{\alpha}\\
			\text{arg}_u(x) + \frac{\tau}{\alpha} &\mbox{if } \quad\text{arg}_u(x)-f(x) < 
			-\frac{\tau}{\alpha}\\
			f(x) &\mbox{if } \quad|\text{arg}_u(x)-f(x)| < \frac{\tau}{\alpha}
			\end{cases}$.
			\item $\bar{u}^{n+1} = u^{n+1} + \theta (u^{n+1} - u^{n})$.
			\item Set $n=n+1$.	
		\end{enumerate}\\
		\State \textbf{end while}\\
		\end{quote}
		\begin{enumerate}\itemsep5pt
		\setcounter{enumi}{1}
		\item Set $u(t_{\alpha_i},x)=u^n$.
		\item Set $\phi(t_{\alpha_i},x) = u(t_{\alpha_{i-1}},x)-u(t_{\alpha_i},x) $.
		\item Set $S(t_{\alpha_i}) = \langle \phi(t_{\alpha_i},x), f(x)\rangle$.
		\end{enumerate}\\
		\State \textbf{end for}\\
		\State \Return $u(t_{\alpha_1},x),...,u(t_{\alpha_N},x)$, $\phi$, $S$.
	\end{algorithmic}
	}
\end{algorithm}
\vspace*{-1.5em}%
The minimization algorithm proposed by Chambolle and Pock consists of three update steps: the first step is a dual update using the resolvent operator of $F^{\ast}$ and the second is a primal update using the resolvent operator of $G$. These are followed by a simple weighting step between the previous two primal iterates. See \cite{Chambolle2011} for more details. The resolvent operators for $G$ and $F^{\ast}$ are presented in \cite[chpt. 6.2]{Chambolle2011}. To construct a forward scale space, we solve the $L^{1}-TV$ denoising model with increasing regularization parameter $t_{\alpha_1} < t_{\alpha_2} < ... < t_{\alpha_N}$ and compute the spectral transform function $\phi(t,x)$ via backward-differences and the response function $S(t)$ based on these solutions. Both the resulting primal-dual algorithm to minimize \eqref{eq:l1tv} and the computation of the spectral functions are embodied in Algorithm \ref{alg:cpforl1tv}. Note that to receive a very high degree of convergence we needed a very large number of iterations resulting in a long computational time. With fewer iterations, the algorithm did not converge completely and eigenfunctions lost contrast altough, according to \eqref{eq:l1sol}, this should not be the case.
%%%%%%%%%%%%%%%%%%%%

\vspace*{-0.5em}%
\section{Results}\label{sec:res}
%%%%%%%%%%%%%%%%%%%%
%\vspace*{-.5em}%
In this section, we describe the main properties and advantages of the $L^1-TV$ based spectral approach introduced above. We discuss both several synthetic experiments and real experiments from biological cell imaging and retina imaging. For some examples, the results are compared to results from the $L^2$ based spectral framework to illustrate the differences between both models. 
For all results, we applied the framework summarized in Algorithm \ref{alg:cpforl1tv} and cluster the spectral response function afterwards. In our experiments, we set $\tau = 0.2, \sigma = 0.625, \theta = 1$ and $maxIts = 50.000$. For synthetic datasets we used $N = 20$ linearly spaced $t_{\alpha_i}$ and increased $N$ to 50 for the experimental datsets. After solving the variational model with increasing regularization parameter, we manually clustered $S$ but comparable results can be achieved with common histogram thresholding methods such as Otsu's method, the Triangle method or methods designed for more than two classes. The different classes in the histogram are always visualized by different colors. Reconstruction of the filtered signal was performed via \eqref{eq:filter} where $H(t)$ was defined as an indicator function of the filtered time interval and the resulting signal was color-coded with the same color as in the spectral response and multiplied with the original gray values so that intensity changes in the input signal remain detectable.\\
\begin{figure}[htb]
  \centering 
    \subfloat[Input Data.]{\includegraphics[height=0.143\textheight]{img/balls_size/orig_data.png}}\quad 
  \subfloat[Spectral Response.]{\includegraphics[height=0.143\textheight]{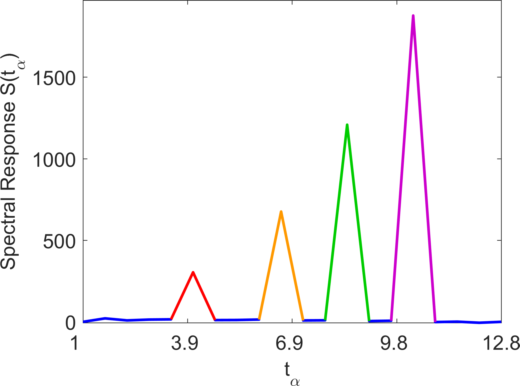}}\quad 
  \subfloat[Color-coded \newline Reconstruction.]{\includegraphics[height=0.143\textheight]{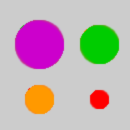}}
  \vspace*{-.5em}%
  \caption{\textit{Detection of size scales for eigenfunctions with constant intensity.} (a) shows the original input data, (b) the resulting spectral response function of the forward $L^{1}-TV$ denoising approach. Every peak in $S$ corresponds to one disc. The color-coded reconstruction is shown in (c). }
\label{fig:ballssize} 
\end{figure}
\begin{figure}[thb]
  \centering 
    \subfloat[Input Data.]{\includegraphics[height=0.143\textheight]{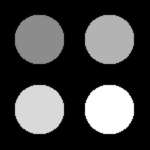}}\quad 
  \subfloat[Spectral Response.]{\includegraphics[height=0.143\textheight]{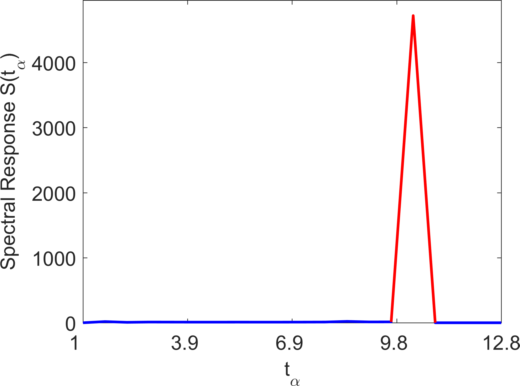}}\quad 
  \subfloat[Color-coded \newline Reconstruction.]{\includegraphics[height=0.143\textheight]{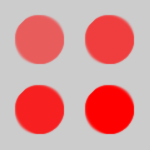}}
  \vspace*{-.5em}%
  \caption{\textit{Invariance w.r.t. to intensity scales.} (a) shows the input data consisting of four TV eigenfunctions of similar size but various intensities. The corresponding spectral response is shown in (b). All discs vanish at the same moment, independent of their contrast. (c) The reconstruction of the peak is shown in red.}
\label{fig:ballsint}%
% \vspace*{-2em}
\end{figure}%

\textit{Detection of Varying Eigenshapes: Size versus intensity.} In this set of experiments, we focus on circular objects that are all eigenfunction of the total variation functional. The aim of these experiments is to investigate which scales can be detected and reconstructed using our method and where the differences to an $L^2$ based approach are. Figure \ref{fig:ballsint} shows an example with four discs in front of a uniform background, all showing the same intensity level. The corresponding spectral response function (b) clearly shows the four peaks each corresponding to one disc, as can be seen in the color-coded reconstruction in (c). Since $L^1$ is a purely size-based approach, it decomposes these eigenfunctions clearly. The opposite case is shown in Figure \ref{fig:ballsint}. Here, the input data is again composed of four discs but now with only one size and various intensity levels. Although an $L^2$ based approach would again show four peaks in the spectral response function due to its size/intensity ambiguity mentioned previously, we now see only one peak in the response (b). Filtering only the signal belonging to this red peak returns all four discs; see (c). This example clearly shows the contrast invariance of the $L^1$ data fidelity term which is a major difference between $p=2$ and $p=1$. 
\begin{figure}[htb]
  \centering 
    \subfloat[Input Data.]{\includegraphics[height=0.143\textheight]{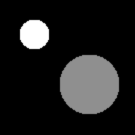}}\quad 
  \subfloat[Spectral Response of \newline $L^2-TV$ denoising.]{\includegraphics[height=0.143\textheight]{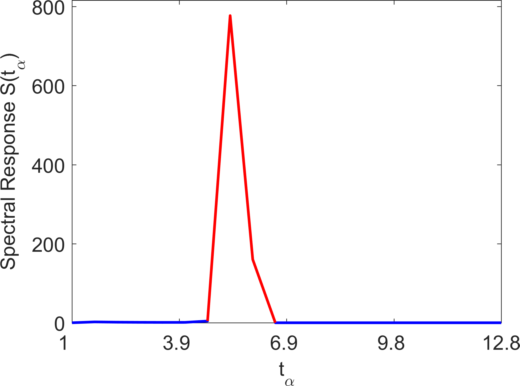}}\quad 
  \subfloat[Color-coded Reconstruction ($L^2$).]{\includegraphics[height=0.143\textheight]{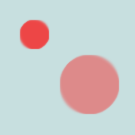}}\\
\subfloat{\includegraphics[height=0.143\textheight]{img/balls_int_and_size/empty.png}}\quad 
\setcounter{subfigure}{3}
  \subfloat[Spectral Response of \newline $L^1-TV$ denoising.]{\includegraphics[height=0.143\textheight]{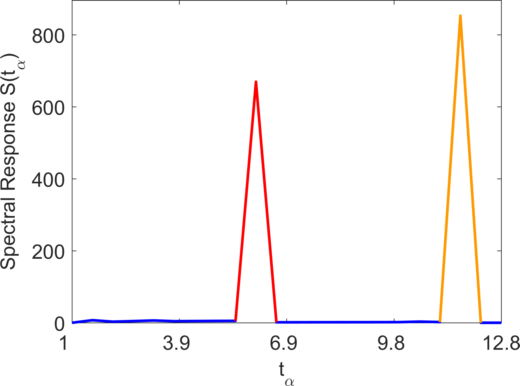}}\quad 
  \subfloat[Color-coded Reconstruction ($L^1$).]{\includegraphics[height=0.143\textheight]{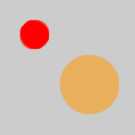}}
  \caption{\textit{Detection of size and intensity scale mixtures.} (a) shows the input data with two eigenfunctions with differing sizes and intensities. Results of an $L^{2}$ based spectral approach are shown in (b) and (c). The method is not able to separate the objects. (d) and (e) show the results of the contrast invariant $L^{1}$ based spectral analysis. The two discs are clearly separable based on size. }
\label{fig:ballssizeandint} 

\end{figure}
\begin{figure}[htb]
% \vspace*{-0.5em}%
  \centering 
    \subfloat[Input Data.]{\includegraphics[height=0.13\textheight]{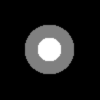}}\quad 
  \subfloat[Spectral Response of $L^2-TV$ denoising.]{\includegraphics[height=0.13\textheight]{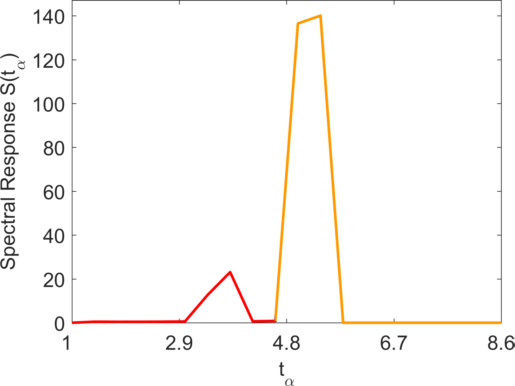}}\quad 
  \subfloat[Color-coded Reconstruction of red peak ($L^2$).]{\includegraphics[height=0.13\textheight]{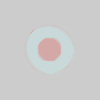}}\quad
  \subfloat[Color-coded Reconstruction of orange peak ($L^2$).]{\includegraphics[height=0.13\textheight]{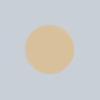}}\\
\subfloat{\includegraphics[height=0.13\textheight]{img/balls_int_and_size/empty.png}}\quad 
\setcounter{subfigure}{4}
  \subfloat[Spectral Response of $L^1-TV$ denoising.]{\includegraphics[height=0.13\textheight]{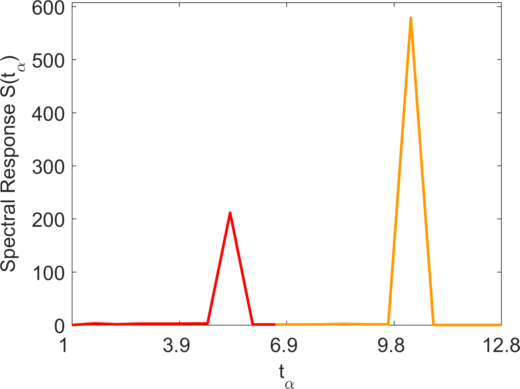}}\quad 
  \subfloat[Color-coded Reconstruction of red peak ($L^1$).]{\includegraphics[height=0.13\textheight]{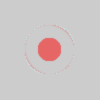}}\quad
 \subfloat[Color-coded Reconstruction of orange peak ($L^1$).]{\includegraphics[height=0.13\textheight]{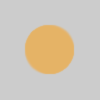}}
  \vspace*{-.5em}%
  \caption{\textit{Overlapping mixtures of size and intensity scales.} (a) shows the input data with two discs with different size and intensity on top of each other. (b) - (d) shows the spectral response and the reconstructions of both peaks using an $L^2$ dataterm and (e) - (g) for an $L^1$ dataterm. In (c) we see artifacts of the larger disc already appearing at the fine scales while (f) and (g) clearly separate the two discs.}
\label{fig:ballsontop} 
 \vspace*{-2em}
\end{figure}
%%%%

A direct comparison of both spectral frameworks is shown in Figures \ref{fig:ballssizeandint} and \ref{fig:ballsontop}. Figure \ref{fig:ballssizeandint} shows as input data (a) two discs with different sizes and intensities. From a visual perspective, this would be clearly identified as two different scales but the $L^2-TV$ denoising approach is not able to distinguish between the objects. The scale that this approach uses is always a mixture of the size scale (small for the disc on top and large for the other one) with the intensity scale (large for the one on top and small for the second disc) and therefore it can occur that they both end up with the same ``medium" scale. In the spectral domain, they are represented by one peak (b) and can therefore not be reconstructed separately. However, this is different for the contrast invariant $L^1-TV$  approach since it is purely size based. The spectral response function shows two easily separable peaks (d) that can be reconstructed one by one (e). 
A similar behavior can be seen for two discs on top of each other (see Figure \ref{fig:ballsontop}(a)). Although these two discs are represented by two separate peaks in both approaches, the $L^2$ based approach mixes them if they are on top of each other. In this case, the peaks become less sparse and apart from each other (b) and in the reconstruction of the small scales (c) we see some artifacts of the larger disc (bluish ring around red disc). This is not the case for $L^1$. Both peaks are clearly separable (e) and the reconstruction gives a clear separation of both discs. The larger red circle in (f) is just a discretization artifact. 

\begin{figure}[htb]
  \centering 
    \subfloat[Input Data.]{\includegraphics[height=0.14\textheight]{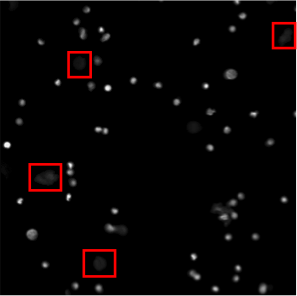}}\quad 
%  \subfloat[Spectral Response.]{\includegraphics[width=0.27\textwidth]{img/cells/S_new_def_2.png}}
  \subfloat[Spectral Response.]{\includegraphics[height=0.125\textheight]{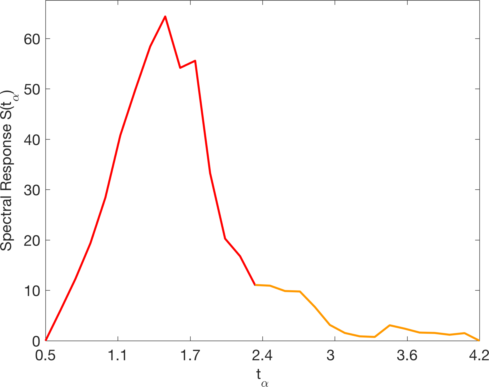}}\quad
  \subfloat[Color-coded Reconstruction.]{\includegraphics[height=0.14\textheight]{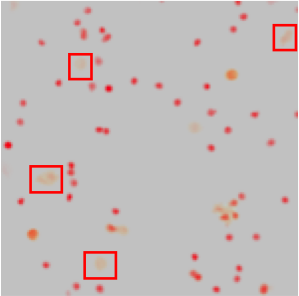}}\quad
  \subfloat[Segmentation using the Bregman-CV.]{\includegraphics[height=0.14\textheight]{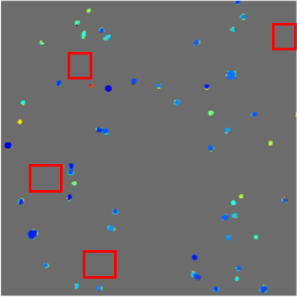}}\\
 \caption{\textit{Segmentation of cells of multiple sizes and intensities.} (a) shows an experimental dataset of tumor cells under a fluorescent microscope; (b) the resulting spectral response function; and (c) the color-coded reconstruction. A segmentation result taken from \cite{Zeune2016} is given in (d). Wee see that both methods can be used to obtain a  multi-scale segmentation but very dim objects (red boxes) might be lost when taking an intensity based approach as in (d). In taking our purely size-based approach the dim cells are also found.}
\label{fig:cells} 
% \vspace*{-0.5em}
\end{figure}

\begin{figure}[htb]%\vspace*{-1.7em}%
  \centering 
  \subfloat{\includegraphics[height=0.12\textheight]{img/balls_int_and_size/empty.png}}\quad 
  \setcounter{subfigure}{0}
  \subfloat[Unbiased Input.]{\includegraphics[height=0.125\textheight]{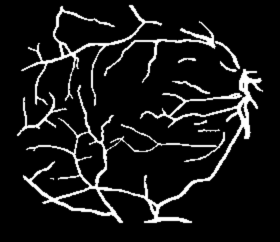}}\quad 
  \subfloat[Biased Input.]{\includegraphics[height=0.125\textheight]{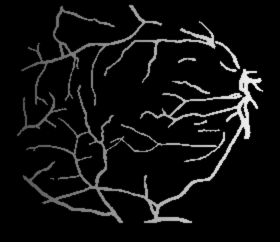}}\quad
  \subfloat{\includegraphics[height=0.12\textheight]{img/balls_int_and_size/empty.png}}\\%
  \vspace*{-1em}%
\setcounter{subfigure}{2}
    \subfloat[Spectral Resp. ($L^2$, Unbiased)]{\includegraphics[height=0.115\textheight]{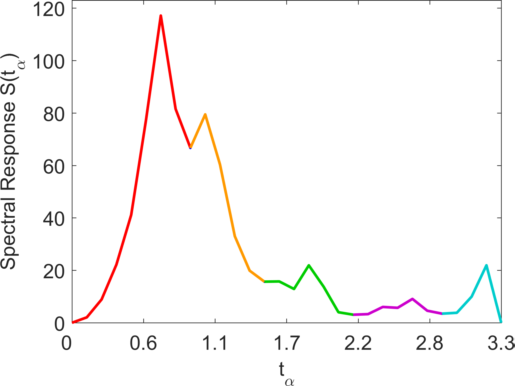}}\quad 
  \subfloat[Color-coded Reconstruction. ($L^2$, Unbiased)]{\includegraphics[height=0.115\textheight]{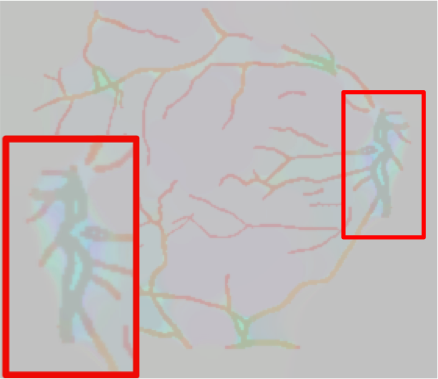}}\quad 
  \subfloat[Spectral Resp. ($L^2$, Biased)]{\includegraphics[height=0.115\textheight]{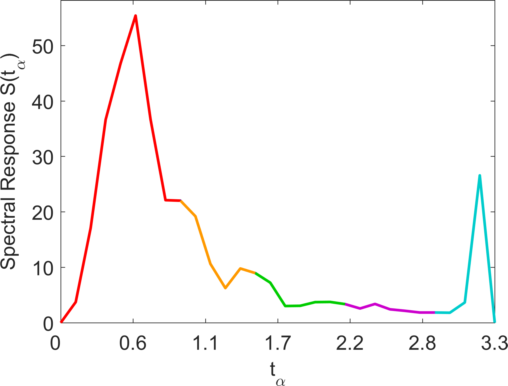}}\quad
  \subfloat[Color-coded Reconstruction. ($L^2$, Biased)]{\includegraphics[height=0.115\textheight]{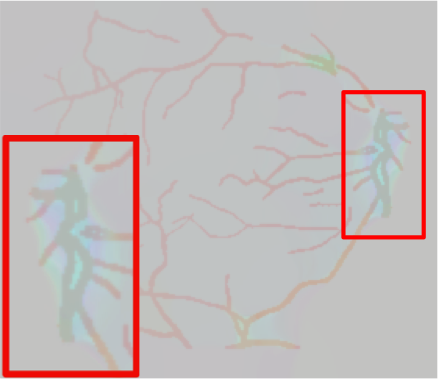}}\\%
  \vspace*{-1em}%
    \subfloat[Spectral Resp. ($L^1$, Unbiased)]{\includegraphics[height=0.115\textheight]{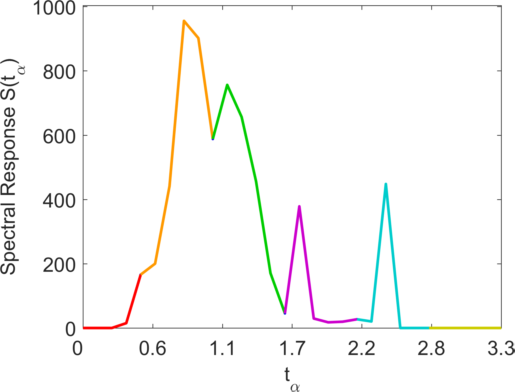}}\quad 
  \subfloat[Color-coded Reconstruction. ($L^1$, Unbiased)]{\includegraphics[height=0.115\textheight]{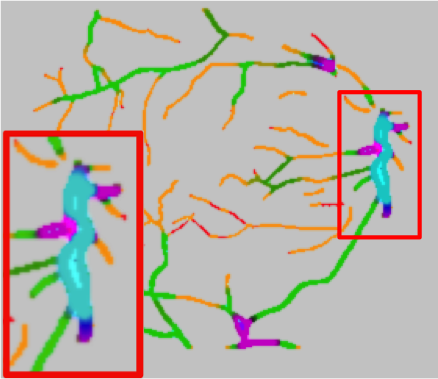}}\quad 
  \subfloat[Spectral Resp. ($L^1$, Biased)]{\includegraphics[height=0.115\textheight]{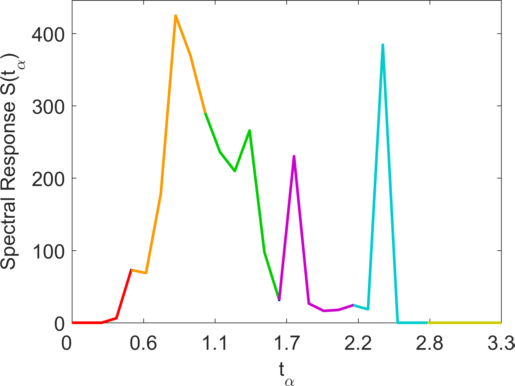}}\quad
  \subfloat[Color-coded Reconstruction. ($L^1$, Biased)]{\includegraphics[height=0.115\textheight]{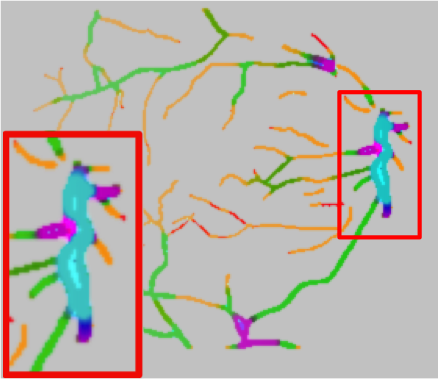}}
  \caption{\textit{Decomposition of network structures without and with intensity bias.} (a) shows the binary input network of blood vessel \cite{hoover2000} and (b) the same network with an intensity bias. (c)-(f) show spectral response and reconstruction results using $L^2$ and the unbiased (c)-(d) resp. biased data (e)-(f). In (g)-(j), the results with $L^1$ are presented. Reconstructions (d) and (f) appear very dim since vessels are not removed in one step but are reshaped over time, resulting in a mix of (colored) scales. This is more clear in the magnification. The turquoise part in (d) and (f) is no longer vessel-shaped but more a roundish shadow around the original shape while the structure remains unchanged for $L^1$ based reconstruction (h) and (j). This network structure is also not affected by an intensity bias.}
\label{fig:network} 
 \vspace*{-2em}
\end{figure}

%%%%
\textit{Segmentation of Experimental Cell Data.} In Figure \ref{fig:cells}, we present a dataset that was also used in \cite{Zeune2016}. The experimental dataset (a) shows a fluorescent microscopy image of Circulating Tumor Cells. Here, the goal is to reliably segment all cells although they differ much in size and intensity. A multi-scale segmentation approach was presented in \cite{Zeune2016}, see (d), but due to the intensity dependency of this approach, the method was not able to detect those cells that are very dim (highlighted with red boxes). Since our new $L^1$ based spectral approach ignores intensity differences, the method also finds the very dim cells. When reconstructing the orange part of the spectral response function (b) and using the thresholding formula \eqref{eq:seg}, we obtain a segmentation that contains all four cells highlighted with a red box.
\newpage
\vspace{-3em}%
\textit{Experimental Data of Network Structures.}  In Figure \ref{fig:network} (a), a manually segmented blood vessel network taken from the STARE dataset \cite{hoover2000} is shown. A problem often occurring in retinal blood vessel segmentation is that small vessels are also very dim and therefore even more challenging to detect. In (b), we added an intensity bias to the original data to test whether this influences our segmentation/clustering approach or not. In (c)-(f), the spectral response functions and reconstructions for $L^2$ denoising are shown. We see that there are fewer clear peaks in the spectral response function, especially in the case of intensity biased input. Another problem that we observe is that the vessels are not removed while retaining their original shape but are reshaped to more circular objects. This leads to a mixing of all scales in the reconstruction and the dim appearance of the visualization. For $L^1$ based denoising, the spectral responses are much sparser in (h) and (j) and therefore easier to cluster. The reconstructions (g) and (j) are also much easier to interpret since the vessels are not reshaped but removed in one step based on their diameter. Thus, for both networks (even for the one with an intensity bias) we see a clear clustering of blood vessels based on size. The difference between both methods is clearly shown in the magnified box. While in the top row the largest scale is just a turquoise shadow around the thickest blood vessels, the blood vessel is clearly reconstructed in the bottom row.
%%%%%%%%%%%%%%%%%%%%
 \vspace*{-0.5em}
\section{Conclusion and Outlook}
 \vspace*{-0.5em}
%%%%%%%%%%%%%%%%%%%%
In this paper, we have described the study of contrast-invariant $L^1$ data fidelities for variational multi-scale methods in combination with nonlinear spectral image analysis. We have shown that the contrast invariance results in an improved sparsity of spectral responses. In comparison to standard $L^2$, this allows a more informative spectral image representation to be obtained. We presented a model, an efficient algorithm and numerical results. In the particular case of experimental data sets that have complex shapes and strong intensity variations of objects or the background, our method outperforms the current standard method for nonlinear spectral decomposition. For future studies, it will be important to extend the ideas to nonlocal graph-based problems such as $L^1$ with nonlocal $TV$, and further analyze the relationship of the doubly nonlinear scale-space flow to the proposed scale-space procedure.
%%%%%%%%%%%%%%%%%%%%
 \vspace*{-0.5em}
\section*{Acknowledgements}%
 \vspace*{-0.5em}
LZ, LT, CB acknowledge support by the EUFP7  program \#305341 CTCTrap and the IMI EU program \#115749 CANCER-ID. CB acknowledges support for his TT position by the NWO via Veni grant 613.009.032 within the NDNS+ cluster.
\vspace*{-1em}%
\bibliographystyle{splncs}
\bibliography{refs}
\vspace*{-1em}%
\end{document}

%% file: commands.tex
\newcommand{\eqn}[1]{\begin{equation}#1\end{equation}}

\newcommand{\R}{\mathbb{R}}

\DeclareMathOperator*{\esssup}{\textnormal{ess~sup}}

\newcommand{\norm}[1]{\left\| #1 \right\|}

\renewcommand{\div}{\nabla \cdot}%
\newcommand{\grad}{\nabla}%